%% file: PredictingCardinalityOfGB-arxiv2.tex
\documentclass{amsart}

\usepackage{amsmath, amssymb,amsthm,amscd,amsfonts,graphicx}%
\usepackage{fancyhdr, geometry}
\usepackage{url, hyperref, color, verbatim,wrapfig}
\usepackage{mathtools,tikz,caption,subcaption}
\usetikzlibrary{matrix}
\usepackage{url, hyperref, color, verbatim}
\usepackage[normalem]{ulem}
\usepackage{natbib}

\theoremstyle{plain} 

\theoremstyle{definition}

\usepackage{tikz}
\usepackage{stmaryrd}

\tikzstyle{dense} = [rectangle, draw, fill=blue!20,
 text width=14em, text centered, rounded corners, minimum height=1.5em]
\tikzstyle{input} = [rectangle, draw, fill=green!20,
 text width=8em, text centered, rounded corners, minimum height=1.5em]
 \tikzstyle{conv} = [rectangle, draw, fill=red!20,
  text width=12em, text centered, rounded corners, minimum height=1.5em]

\oddsidemargin \evensidemargin
\newcommand{\ideal}[1]{\left< {#1}\right>}

\DeclareMathOperator{\init}{LT_{\prec}}

\def\jdlqed{\vbox{\hrule \hbox{\vrule\hbox to
5pt{\vbox to 6pt{\vfil}\hfil}\vrule}\hrule}}

\title{\bf Predicting the cardinality and maximum degree of a reduced Gr\"obner basis}
\author{Shahrzad Jamshidi, Eric Kang, and Sonja Petrovi\'c}

\date{}

\begin{document}

\begin{abstract}
We construct neural network regression models to predict key metrics of complexity for Gr\"obner bases of binomial ideals. This work illustrates why predictions with neural networks from Gr\"obner computations are not a straightforward process. Using two probabilistic models for random binomial ideals, we generate and make available a large data set that is able to capture sufficient variability in Gr\"obner complexity. We use this data to train neural networks and predict the cardinality of a reduced Gr\"obner basis and the maximum total degree of its elements. While the cardinality prediction problem is unlike classical problems tackled by machine learning, our simulations show that neural networks, providing performance statistics such as $r^2 = 0.401$, outperform naive guess or multiple regression models with $r^2 = 0.180$.
\end{abstract}


\maketitle

\section{The problem}

Let  $K$ be a field---for example  the reader may keep $K=\mathbb C$ in mind.  Let $f_1=0,\dots,f_s=0$ be a system of $s$ polynomial equations in $n$ variables with coefficients in $K$.  
We will denote by  $\ideal{f_1,\dots,f_s}\subset R=K[x_1,\dots,x_n]$ the ideal generated by these polynomials. If $F=\{f_1,\dots,f_s\}$ is a set of polynomials, the ideal $\ideal{f_1,\dots,f_s}$ will equivalently be denoted by $\ideal{F}$. 
 Fix a monomial order $\prec$ on $R$ and an ideal $I$.  
Denote by $\init(h)$ the initial  term  of $h\in R$, which is the $\prec$-largest monomial. The initial ideal of $I$ with respect to $\prec$ is the monomial ideal $\init(I):=\ideal{\init(h):h\in I}$. 

A \emph{Gr\"obner basis} of $I$ with respect to $\prec$ is a set $G=\{g_1,\dots,g_k\}$ such that  $\prec$-largest monomials $\init(f)$ satisfy  $\init(f) \in \ideal{\init(G)}$ for all $f\in I$.  While a Gr\"obner basis, or its size, is not unique for a given order,  a reduced one---where redundant terms are removed from each $g_i$---is.
 \cite{StWhatIsGB} offers a high-level overview of Gr\"obner bases and the textbook \cite{clo} discusses various applications as well. Of the many applications we single out two that have generated tremendous interest in recent decades: discrete optimization \cite{GBandOptim,rekha} and statistics \cite{GBandStats, DS98}.

The starting point of our work are the following  questions:
\begin{enumerate}
  \item[1: ] \label{problem1}
   {\it For a given monomial ordering, can we reliably predict the size of the reduced
  Gr\"obner basis of an ideal?}
  \item[2: ] \label{problem2} {\it For a given monomial ordering, can we reliably predict the maximum total degree
      of the reduced Gr\"obner basis of an ideal?
      }
\end{enumerate}
A key word here is \emph{predict}; we show as a proof of concept  that the size of a Gr\"obner basis of an ideal can, in fact, be predicted by a machine learning (ML) algorithm.
 \cite{DylanMike-LearningBuch}  and \cite{MPP22} show how
 a handful of cornerstone algorithms in computational algebra are well-suited
 for a machine learning approach. Our work follows this promising direction of exploration. 
 While the use of machine learning is under-explored in the particular computations  in which we are interested, there is a wider literature on the use of machine learning for algebraic problems. For example,  \cite{MLMathStructures} offers an overview of recent research on machine-learning mathematical structures; see also \cite{DeepLearnForSymbolicMath,DeepLearnForMath}. The thesis \cite{Silverstein}  uses features in neural network training to select the best algorithm to perform a Hilbert series computation by predicting a best choice of a monomial for one of the crucial steps of the computation, namely predicting a best pivot rule for the given input.

In this paper, we answer the prediction problem 1  affirmatively for the case of binomial ideals  in $5$ variables under graded reverse lexicographic  order. We also attempt to answer problem 2, but the performance is less good, so we suggest possible avenues for further exploration. 
We will quantify what `reliable' prediction means, describe the data generating process, and  compare the machine learning model performance with other methods. Along the way, we  discuss further potential developments on the computational front.

\subsection{Problem relevance}

One  problem whose solution critically depends on computing a Gr\"obner basis is the ideal membership problem (IMP),
which asks:
\emph{Given the polynomials $\{f, f_1, \ldots, f_s\} \subset k[x_1, \ldots, x_n]$, is $f\in  (f_1, \ldots, f_s)$?}
A surprising number of applied computational questions can reduce to the ideal membership problem (IMP). For example, the question, ``Does a system of polynomial equations $\{f_1, \ldots, f_s\}$ have a set of common roots?'' can be rewritten as an IMP: is $1 \in ( f_1, \ldots, f_s )$?
Once a Gr\"obner basis of  the ideal $( f_1, \ldots, f_s )$ is computing, the IMP can be solved using polynomial division.
 Many computer algebra systems offer generic algorithms for computing Gr\"obner bases applicable to all kinds of input ideals. Despite this great success in the field, unfortunately, algebraic problems have well-known bad worst-case complexity.   In 1965, \citet{Buchberger1965Translated} proposed  a groundbreaking algorithm which is able to compute a Gr\"obner basis of {\it any} ideal. Unsurprisingly, as the problem is $NP$-hard in general, it has a doubly exponential expected runtime in the number of variables~\citep{dube}, although
 in a few cases, specialized algorithms have been used to improve runtime~\citep{BePa08, BePa08a}.
 Challenges in computing Gr\"ober bases give rise to an interest in exploring how new techniques might provide alternatives.
 Many algorithmic tasks  in symbolic computation relate to solving polynomial systems, a task which also requires a computation of a Gr\"obner basis; this 
 relies upon variants of Buchberger's approximately 60-year-old algorithm.

Recently, \cite{GBviolator} proposed a new machine-learned randomized framework, called Spark Randomizer, for computing Gro\"bner bases. The approach is distinct from traditional methods because it relies on a cleverly biased random sampling method rather than cleverly organized versions of Buchberger’s algorithm. This framework takes a departure from the standard algorithms based on S-polynomials and address the problem of computing Gr\"obner bases using violator spaces, a concept used in geometric optimization. 
One key input to this framework is an estimate of the cardinality of a minimal Gr\"obner basis, that is, the answer to Problem 1. 

Two remarks are in order. First, it seems that, at the moment, one cannot provide predictions with extremely high accuracy. In these cases, the consequences of inaccurate estimation can be two-sided. In case of over-estimating, the  Spark Randomizer will still produce correct output, albeit not necessarily minimal. In case of under-estimating, it will simply never produce an output! While this isn't exactly an exciting prospect, at least it does not lead to incorrect answers in terms of Gr\"obner bases.
Moreover, purely as a theoretical exercise, predictions from problems 1 and 2 can be useful for, for example, informing which method or monomial ordering to use for computing a Gr\"obner basis.

\subsection{Difficulty assessment: an ML perspective}
\begin{figure}[h]
    \includegraphics[scale = 0.22]{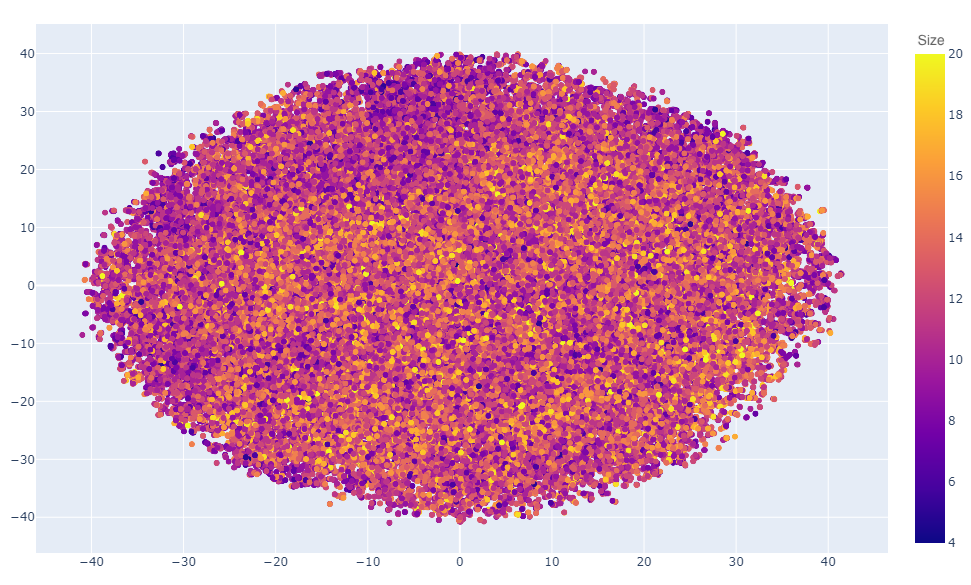}
 \caption{ A t-SNE plot for $99,622$ binomial ideals in $3$ variables, each generated by $5$ homogenous binomials with total degree $7$. Output color represents the size of the reduced graded reverse-lexicographic Gr\"ober basis.}
				  \label{fig:tSNEhom}
\end{figure}
Justifiably, one might wonder why the problems we are posing are not just simple tasks that can be completed using out-of-the-box machine learning tools.  In other words,
why is this prediction problem not of `plug-and-play' type?
The difficulty lies in the fact that the Gr\"obner size or degree bound problems are very much \emph{unlike} traditional machine learning benchmark problems, in which the input is complex but the humans can quickly validate the prediction.

As an illustration, consider a data set of $99,622$ binomial ideals in $3$ variables, each generated by $5$ homogenous binomials with total degree $7$. For a human to `quickly validate the prediction' one either needs good folklore knowledge of the problem, which does not exist for a random ideal, or one needs a way to visualize this high-dimensional data in two or three dimensions. A common statistical method for such a visualization is called t-distributed stochastic neighbor embedding, or a t-SNE plot. It gives each datapoint a location on a 2- or 3-dimensional map; the data are then color-coded by the output, in our case, size of the reduced graded reverse-lexicographic  Gr\"obner basis.
We represent each ideal through the natural exponent vector representation  of the generators described in Section \ref{sec:data}.
The t-SNE plot  of this data set in Figure~\ref{fig:tSNEhom}  demonstrates a lack of clustering and the absence of a decision boundary, suggesting a difficult classification problem.

A practitioner in machine learning who understands Gr\"obner basis size does not find the t-SNE plot surprising.
Namely, the t-SNE algorithm relies on the Euclidean distance between these vectors and converting them to conditional probabilities in order to visualize data. This is not an arbitrary choice; this data visualization technique is relevant for assessing problem difficulty from the point of view of ML precisely because  the $L^2$ norm, or mean squared error, is commonly used in training neural networks via back propagation (i.e., gradient descent).
However, this notion of distance is meaningless for these vectors.
For example, if we consider three particular ideals from the t-SNE plot above, $I$, $J$ and $K$, whose 
reduced graded reverse lexicographic Gr\"obner bases are of sizes 7, 13, and 14, respectively.
Their vector representations are explicitly stated in Table~\ref{tab:distanceNonsenseExample}.
\begin{table*}[h]
\begin{tabular}{rl}
    Ideal: & $I=\ideal{x^3y^3z-xy^3z^3, x^6y-x^2yz^4, xy^5z-x^4yz^2,y^6z-x^3y^2z^2,x^6y-xy^3z^3}$ \\
   Vector representation: &     $v_I=(3, 3, 1, 1, 3, 3, \phantom{x}  6, 1, 0, 2, 1, 4, \phantom{x} 1, 5, 1, 4, 1, 2, \phantom{x} 0, 6, 1,3, 2, 2,\phantom{x}  6, 1, 0, 1, 3, 3)$ \\
   Size of reduced GB(I):     & 7  \\
   \hline
    Ideal: & $J = \ideal{x^6z-x^3y^2z^2, \phantom{x} xy^3z^3-xz^6,\phantom{x}  y^7-y^4z^3,\phantom{x} x^4yz^2-xyz^5,\phantom{x} xy^6-y^3z^4}$ \\
   Vector representation: &     $v_J=(6, 0, 1, 3, 2, 2, \phantom{x} 1, 3, 3, 1, 0, 6, \phantom{x} 0, 7, 0, 0, 4, 3, \phantom{x} 4, 1, 2, 1, 1, 5, \phantom{x} 1, 6, 0, 0, 3, 4)$ \\
   Size of reduced GB(J):     & 13  \\
   \hline
    Ideal: & $K = \ideal{x^5y^2-yz^6, \phantom{x} x^6y-x^2y^2z^3, \phantom{x} x^3y^3-x^1yz^3,\phantom{x} y^3z^4-xz^6,\phantom{x} x^3y^4-y^7}$\\
   Vector representation: &    $ v_K=(5, 2, 0, 0, 1, 6,\phantom{x}  6, 1, 0, 2, 2, 3,\phantom{x}  3, 4, 0, 3, 1, 3, \phantom{x} 0, 3, 4, 1, 0, 6, \phantom{x} 3, 4, 0, 0, 7, 0)$\\
   Size of reduced GB(K):     & 14 \\
   \hline
\end{tabular}
\medskip
\caption{Three ideals in $3$ variables, each generated by $5$ homogenous binomials of degree $7$, their vector representations, and the sizes of their reverse-lexicographic reduced Gr\"obner bases, from the data set depicted in the t-SNE plot in  Figure~\ref{fig:tSNEhom}. In the table, `GB' stands for Gr\"obner basis.}
\label{tab:distanceNonsenseExample}
\end{table*}
It so happens that the Euclidean distance between their vector representations
$v_I,v_J,v_K$ satisfies $d(v_I, v_K) < d(v_J, v_K) < d(v_I, v_J)$.
This shows that the distances do not relate to the GB sizes. 

\subsection{A family of specific problem instances}
As stated, the general Gr\"obner size prediction problem is too broad, so specifying a particular family of ideals to study is a necessity.
In what follows, we will answer the prediction problems 1 and 2  affirmatively for the case of binomial ideals  in $5$ variables under graded reverse lexicographic  order.

The case for choosing the class of  binomial ideals  has been well argued  in commutative algebra, as this class embodies all the richness and complexities of computations with ideals. A recent summary of suitability of binomial ideals for learning can be found in \cite{DylanMike-LearningBuch} from the point of view of reinforcement learning: some of the hardest polynomial problems are binomial, and they can be generated randomly to avoid generic, or uninteresting, behavior, which in this instance means avoiding zero-dimensional ideals where one expects Gr\"obner bases to be easier. \cite[Section 2.2, Figures 1, 2, and 3]{MPP22} further illustrates that binomial ideals accurately capture much of the Gr\"obner basis problem: there is a large variance in difficulty within distributions,  difficulty increases as expected when we increase the number of variables, and mostly as expected when we increase the number of generators.

In general, we work in the space of binomial ideals in polynomial rings with few variables, between $3$ and $7$. While we have performed the learning problem on polynomial rings between $3$ and $7$ variables, in this paper we focus on $5$ variables to keep as much of the complexity as we are able while being able to run millions of Gr\"obner computations on a desktop computer.

In terms of machine learning,
fixing the number of variables and maximum total degree, along with restricting all coefficients to units, provides the following advantages:
\begin{itemize}
  \item There is a precedent for similar restrictions in other Gr\"obner learning problems \cite{DylanMike-LearningBuch}, \cite{DylanPhD},  \cite{MPP22}, all of which have used $n=3$ variables; see also \cite{Silverstein} and \cite{MLMathStructures};
  \item We are able to  compute millions of reduced Gr\"obner bases for such cases
  on a home computer using {\tt Macaulay2}; and
  \item There is a straightforward and lossless correspondence with degree vectors or matrices
  because the coefficients are restricted to $\pm1$.
\end{itemize}
This last point is the most compelling, 
 as it  implies that  the information necessary to compute the combinatorial dimension
is preserved when representing a binomial ideal as a vector or matrix of generator exponents.
This ensures an invariant exists within the data. It would be possible to expand this work, however, to include coefficients beyond $\pm1$, but it would add some complexity to the learning problem---assuming, of course, the neural network is doing a computation in some way comparable to computing a Gr\"obner basis.

The choice of graded reverse lexicographic monomial order follows the default choice  in {\tt Macaulay2}, which we used to generate training data.

\section{Data and Representation}\label{sec:data}

On a first glance, producing a reasonable training set may seem quite challenging.
We rely on a random binomial  model similar to the ``n-d-s-uniform" model from \cite{MPP22}.
 This random model is built using three basic parameters associated with a polynomial ideal that closely related to computational difficulty: the number of variables $n$,  the maximum total degree $d$ of a monomial in the support of the generators, and the number of  generators $s$.
To generate one ideal in the sample, the model randomly constructs $s$ binomials, each by sampling pairs of two distinct  monomials from the set of all nonzero monomials in $n$ variables up to (or exactly equal to) degree $d$.  This  model uses insights from random monomial ideals \citep{rmi} and  has been shown in \cite{MPP22} to generate data that is able to capture sufficient variability in the Gr\"obner complexity problem.

We generate three data sets as follows.
 In data set A, monomials are selected uniformly at random from the set of monomials of total degree \emph{exactly} $15$. In data set B, monomials are selected uniformly at random from the set of monomials of total degree \emph{less than or equal to} $15$. Since there are many more monomials to sample from and we use simple data generating mechanisms, this data set is smaller than A due to the code scalability paywall. Data set C is a representation of  B, in which the binomial generators themselves are replaced by the following  summary statistics: minimum, maximum, mean, and variance of generator degrees, total number of generators, and the  dimension and degree of the corresponding variety.
\begin{table}[h]
  \centering
  \begin{tabular}{c|cc}
    \hline
    Data set & Model & Sample size \\
    \hline
    A & 5-5-15-homogeneous & 1,000,000\\
    B & 5-5-15-general & 500,000\\
    C & 5-5-15-features & 500,000\\
    \hline
  \end{tabular}
  \caption{Data sets of binomial ideals. Each data set contains ideals generated by $5$ binomials in $5$ variables of degree at most $15$. Data set C is a feature representation of data set B.}
  \label{table:data}
\end{table}
With data set C, we seek to compare the methodology of directly representing the data to the network vs the methodology of doing so via human engineered features.  Simple feature learning  is in line with the results in \cite{Silverstein}, which saw success utilizing features in neural network training to select the best algorithm to perform a Hilbert series computation by predicting a best choice of a monomial for one of the crucial steps of the computation, namely predicting a best pivot rule for the given input.

The information we wish to predict---the size and the maximum total degree of the reduced Gr\"obner basis---can be calculated from the initial generating set of the ideal. Therefore, a unique, lossless encoding of the initial generators ensures the possibility of predicting the target information---the size of the minimal Gr\"obner basis. The encoding used is the degree vector, with the requirement that the binomials are sorted according to term order, which ensures unique data representation. 
 This is similar to flattening the encoding used in \cite{DylanPhD}.

For example, suppose the starting generating set is
\begin{align*}I=(&x^7z^8 - v^2x^2yz^2, \\ &x^{13}y^2  - v^3w^3x^5z, \\ &v^6w^4xy^3 - w^3x^2y^4z^3, \\ &w^2x^5y^2z^2 - wx^2yz^6, \\ &v^3w^2y^5z - v^3wxy^3z^2).\end{align*}
This set of five binomials corresponds to either a 25-dimensional vector or a $10 \times 5$ matrix:
\begin{equation} \begin{pmatrix}\label{eq:representation} 0 & 0 & 7 & 0 & 8 & 2 & 0 & 2 & 1& 2\\ 0 & 0 & 13 & 2 & 0 & 3 & 3 & 5 & 0 & 1 \\ 6 & 4 & 1& 3 & 0 & 0 & 3 & 2 & 4 & 3 \\ 0 & 2 & 5 & 2 & 2 & 0 & 1 & 2 & 1 & 6 \\ 3 & 2 & 0 & 5 & 1 & 3 & 1 & 1 & 3 & 2
\end{pmatrix} \end{equation}
or the equivalent matrix constructed based on the line breaks of the above vector.

The two-dimensional analog (the matrix representation), identical to the encoding used in \cite{DylanPhD}, was tested with no noticeable improvement in the neural network model. Similarly, a three-dimensional analog, with one dimension referencing the polynomial, one dimension referencing the variables, and one dimension referencing the term was also tested. No noticeable improvement was found in this case either.

The reduced reverse lexicographic Gr\"obner basis  of the above example $I\subset \mathbb Q[ v, w ,x ,y ,z]$ has $226$ binomials, and the largest degree of its elements is $29$.
The machine learning algorithm is trained to learn, from a large data set of known Gr\"obner sizes and maximum degrees, to predict these two numerical values from the input vector representation of $I$.

\section{Neural networks} 
In order to understand the process and difficulties in training a machine learning algorithm to answer questions 1 and 2, we offer a brief overview of neural network terminology.

\subsection{Overview}
\input{sectionBackgroundNN.tex}

\subsection{Preliminary Analysis}

Initial investigations were done with data sets with 5 binomial generators of the ideal. The monomials were selected uniformly at random from a set of monomials in 3 variables with total degree {\it exactly} 7 (homogeneous) or {\it at most} 7 (nonhomogeneous). Hence, they were simpler than data sets A, B, and C, mentioned in Table \ref{table:data}. One important difference is that these datasets did not fix a monomial order for the representation. 

Work with these initial datasets provided evidence that the invariant---pattern connecting the input and labels---of the data was highly complex. To what extent is
this invariant discoverable with standard techniques? Top performances with a
testing set size are summarized in the chart in Figure \ref{chart}. We included a
measurement of the percentage of predictions that overshot (guessed a larger
combinatorial dimension) than the true value because this is an acceptable
error for the proposed application. It is important to note that in none of
the techniques used did we develop a specialized measure of loss (when
applicable) to account for this preference.

\begin{wrapfigure}{r}{0.47\textwidth}\label{H_naiveBayes}
    \includegraphics[scale = 0.38]{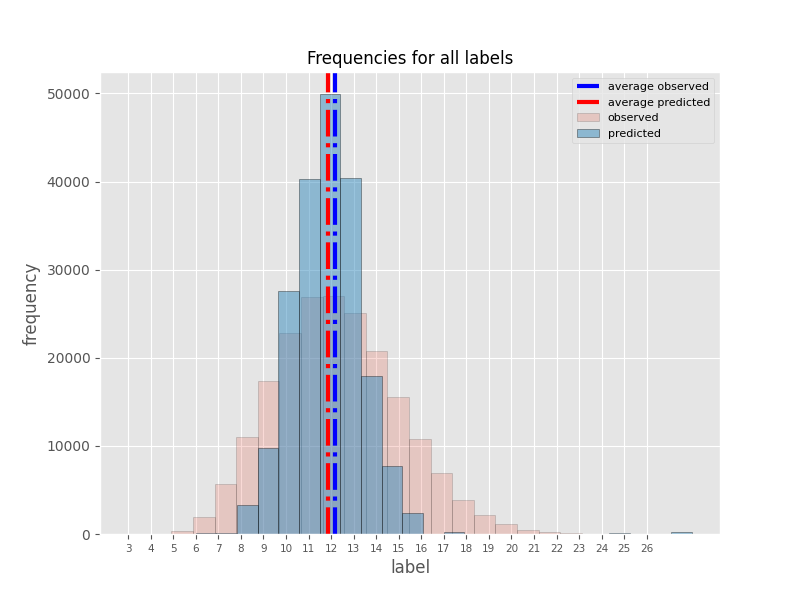}

    \emph{\small A histogram of predicted labels (blue) from the Naive Bayes model compared to the actual labels (pink) of 800,000 samples of homogeneous ideals.}
     \vspace{5mm}
\end{wrapfigure}

All methods performed better than random guessing. We obtained diagnostic plots of performance for each method; but in the interest of space we include only a few that best demonstrate the ideas and issues.

\begin{figure}
\tiny
 \begin{center}
  \begin{tabular}{c||cccc}
    \hline
    \multicolumn{5}{c}{\bf Five Homogenous Binomials in Three Variables with
     Total Degree 7}\\
    \hline
    {\bf Model} & {\bf Training Size} & {\bf \% Accuracy (A)} & {\bf \%
    Overshooting (O)} & {\bf Random A/O}\\
    \hline
    \hline
    {\bf Naive Bayes (multinomial)} & 800,000 & 14.32\% & 53\% & 4.17\%/62\% \\
    {\bf NN} & 796,975 & 14.04\% & 59\% & 5.55\%/52\% \\
    {\bf SVM (RBF)} & 996 & 13.89\% & 51\% & 5.55\%/49\% \\
    {\bf Random Forest} & 49,811 & 8.72\%  & 59\% & 5.88\%/49\%\\
    \hline
  \end{tabular}

  \begin{tabular}{c||cccc}
    \hline
    \multicolumn{5}{c}{\bf Five Nonhomogenous Binomials in Three Variables
    with Total Degree 7}\\
    \hline
    {\bf Model} & {\bf Training Size} & {\bf \% Accuracy (A)} & {\bf \%
    Overshooting (O)} & {\bf Random A/O}\\
    \hline
    \hline
    {\bf Naive Bayes (multinomial)} & 800,000 & 18.82\% & 44\% & 4\%/76\% \\
    {\bf NN} & 798,183 & 24.3\% & 60\% & 6.67\%/64\% \\
    {\bf SVM (RBF)}  &  1,000 & 16.67\% & 50\% &  4.3\%/78\% \\
    {\bf Random Forest} & 49,886 & 10.06\% & 48\% & 7.14\%/64\% \\ \hline
  \end{tabular}
 \end{center}
 \caption{\small Top performing models (by accuracy) based on predictions of a
 testing set. Also included is the percent of predictions that overshot
 (guessed a larger category) as this is an acceptable error. These are compared
 to random guessing (accuracy/overshooting).} \label{chart}
\end{figure}
 \normalsize

\subsubsection{Naive Bayes Classifier}
The multinomial Naive Bayes classifier is suitable for
classification with discrete, integer features, like our data sets. For the
homogeneous data set, we achieved a 14.04\% accuracy. The distribution of
predicted labels, as seen in the figure on the previous page 
follow a normal
distribution much more closely than the actual labels, according to the
Anderson-Darling test.

The analogous nonhomogenous case had similar results. The measured accuracy was
18.82\%, an improvement over the homogeneous case. The only other method to outperform it for this case was the feed-forward neural network. 

\subsubsection{Feed Forward Neural Network}

The first class of neural network methods attempted were LSTMs, a type of
recurrent neural network that is well suited to integral vectors where order
matters. These have had promising results in other algebraic applications
\cite{Silverstein}.
Unfortunately, with these methods we consistently converged to the
most common label or alternating between the two most common (if their
quantities were comparable). Implementing weighting and utilizing various
kinds of optimizers led to the LSTMs to converge to a random label, suggesting that the model was not able to find the invariant.
Networks with more neurons and/or labels proved to be too computationally
expensive for the equipment currently used. So this approach was set aside.  

Feed-forward neural networks (FFNN) performed relatively well compared to other
models investigated. 763 FFNN models were explored in this preliminary investigation with respect to the
homogeneous data sets; from those cases, we were able to
make some valuable observations regarding the number layers, the number of neurons per layer, and optimizers. The best model presented in the chart for both the homogeneous and nonhomogeneous data
consisted of six hidden layers, each with 100 neurons using ReLU activation
(except for the final layer, which used softmax). But our comparisons showed
that networks with 3, 4 or 5 hidden layers performed at comparable levels of
accuracy. The choice of acivation--ReLU with softmax in the final layer, ReLU
with sigmoid in the final layer, sigmoid throughout--also had comparable
performances.

Regarding convergence of the network, it appeared that only 15 to 20 Epochs
was sufficient when plotting training versus validation loss.

\subsubsection{SVM}

Because of computational limitations, the training set of an SVM needs to be
relatively small. Given the complexity of the problem, these seem a
sub-optimal choice; however, we have received questions about SVM performance numerous times. With only 1000
datapoints, an SVM with a radial basis function performed rather well. The primary difficulty with this approach is scale. If this approach were pursued elsewhere, more work would be needed to ensure the training set is both small but sufficiently informative for the space of ideals.

\subsubsection{Random Forests}

Given the categorical nature of both the features and the labels, we explored
random forests with a maximum depth of 10, 100 estimators, and a Gini impurity
criterion. The accuracy was quite low.

\subsection{NN choices for the Gr\"obner prediction problems}\label{sec:our NN}

After experimenting with well over 25,000 trials with each data set listed in Table \ref{table:data}, varying the number of hidden layers, number of neurons, etc., the best results were achieved using the neural network model predicting the size of the reduced Gr\"obner basis using the matrix representation of the degree vector with three hidden layers: a 2D convolutional layer of 300 nodes with $2\times 2$ convolutions, and two dense layers of 500 nodes, with dropout implemented in every layer. 
The topology is illustrated in Figure~ \ref{fig:NNtopology}. This topology and dropout threshold proved to work well for predicting both maximum total degree and the size of the reduced Gr\"obner basis. Additional layers showed no improvement. We tested layers with nodes the ranges of 100 to 2000+. The ranges from 300 to 500 nodes were roughly similar, with a slight improvement observed at 500. The network was trained on 80\% of the data set, with 20\% of the data set reserved as the testing set for benchmarking.

The loss function used was the log hyperbolic cosine function $\sum \log(\cosh(\hat{y} - y))$, which through experimentation proved to be marginally better. In addition, we used the Adam optimizer.  The model was trained for 100 Epochs using a 10\% validation split on the training set to verify that we were not overtraining. Batch sizes of 128 were observed to give the best convergence.

We found no compelling results when providing summary statistics (minimum, maximum, mean, and variance of generator degrees, total number of generators, and the  dimension and degree of the corresponding variety) on the ideals. The matrix representation of the degree vector appears to be the best input for prediction. Similarly, we found no model that predicted the maximum total degree well either.

\begin{figure}[h]
\begin{tikzpicture}[node distance = 4em, auto,ultra thick]

\node [input, rotate=90, ](I) at (0, 1.5) {Input: 15 nodes};
\node [conv, right of=I, rotate=90](C1) {2x2 2DConv 300 nodes};
\node [dense, right of=C1, rotate=90](D1) {500 Nodes};
\node [dense, right of=D1, rotate=90](D2) {500 nodes};
\node [right of=D2, rotate=90, black, draw, fill=green!20, circle, inner sep = 0.25em](D3) { };

\draw[->, ultra thick] (I)--node[above]{\tiny DO 0.5}(C1);
\draw[->](C1)--node[above]{\tiny DO 0.5}(D1);
\draw[->](D1)--node[above]{\tiny DO 0.5}(D2);
\draw[->] (D2)--node[above]{\tiny DO 0.5}(D3);

\end{tikzpicture}
 \caption{\emph{\small The best predictions for both maximum total degree and reduced Gr\"obner basis size was the above topology. The output layer is a single neuron, the value of either degree or size. A dropout threshold of 0.5 was used; it is abbreviated above as ``DO 0.5.'' Because each layer consists of so many neurons, we elected to illustrate the layers of neurons as rectangles, except for the final layer.}
 }
        \label{fig:NNtopology}
\end{figure}

\section{Prediction results}

\begin{figure}[h]
  \includegraphics[scale = 0.35]{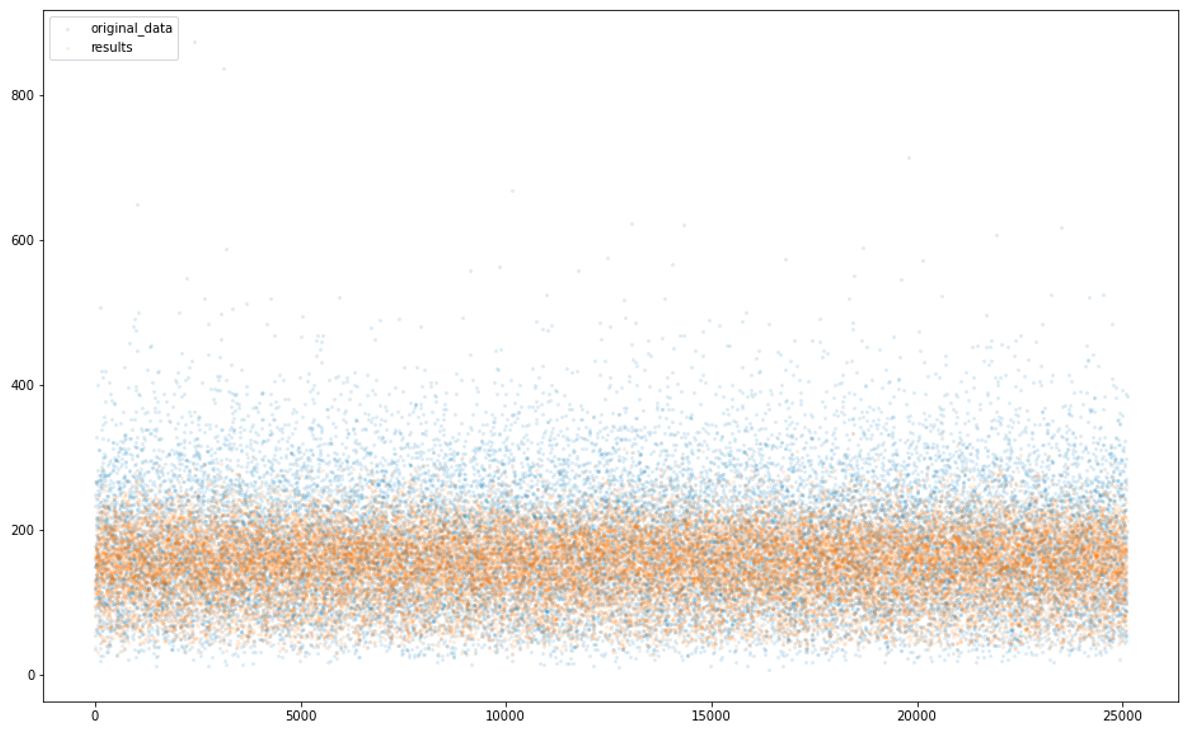}
\caption{Predictions (orange) vs actual values (blue) of the testing set. The $y$-axis represents the size of a reduced revlex Gr\"obner basis for data set B. The $x$-axis represents each trial. Validation r-square is $0.42$. 
} \label{fig:prediction of combinatorial dimension}
\end{figure}

\begin{figure}[h]
  \includegraphics[scale = 0.7]{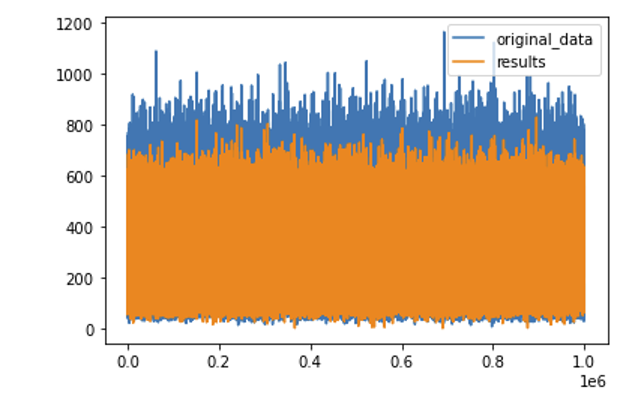}
\caption{Predictions (orange) vs actual values (blue) of the testing set: size of a reduced revlex Gr\"obner basis for data set A. Similar to data set B, the full range of values is not realized in the neural network predictions. Validation r-square is $0.6148$.} \label{fig:prediction of combinatorial dimension for data set A}
\end{figure}

The plots in Figures~\ref{fig:prediction of combinatorial dimension} and \ref{fig:prediction of combinatorial dimension for data set A} shows that the neural network generally did well with more typical sizes of reduced Gr\"obner bases within the data set.  Values on extreme highs and lows were not predicted; the neural network was trying to minimize loss. More complicated techniques are needed to predict these extreme cases.

\medskip
Finally, we offer a summary of performance of the learning algorithms by reporting the $r$-squared statistic in Table~\ref{r2s}.
It is interesting to compare the neural network performance with multiple linear regression models, since that is a reasonable benchmark that can be used to determine if complex learning models (beyond standard regression) are even necessary for the problem at hand.
The $r$-squared statistic, or  the coefficient of determination, is a statistical measure of how close the data lie to the estimated regression line. For example, the baseline model which always predicts the mean will have an $r^2$ value of 0, and the model that exactly matches the observed values will have $r^2= 1$. Models that do worse than then baseline prediction will have a negative $r^2$ value; this often happens when the linear regression is trained (or fitted) on a particular data set and tests or evaluated on a completely different
one (e.g., data draws from a completely different distribution).

The reader may wonder why we attempted to use such simple data features to learn anything about  the Gr\"obner bases. The reason was that a subset of features from data set C  did a reasonably good job in predicting the number of additions during one run of Buchberger's algorithm, which is one of the results in \cite{MPP22}. Thus, from a pragmatic point of view, one could say that there was a precedent--if not hope--that generator features would provide reasonable predictions. They do not, which then in turn is not a surprise to nonlinear algebraists, because performance metrics of Buchberger's algorithm need not correlate with complexity of Gr\"obner bases.

In summary, the neural network model predicted the reduced Gr\"obner basis size better than multilinear regression. It performed especially well when the total maximum degree was fixed at 15. This is, in fact, exactly the outcome one would have hoped for---the learning problem is highly nontrivial, but possible. We found no neural network model that predicted the maximum total degree well nor one that could predict from summary statistics better than multilinear regression. In all cases, there are likely more opportunities for improvement in tackling the learning problem by reconsidering some of the more complex hyperparameters such as the cost function. A cost function that is more meaningful to the problem could allow for significantly better convergence to an optimal outcome.

\begin{table*}[h]
\begin{tabular}{cl|c|c}
Input data & Predicted quantity & Multiple linear regression & Neural network from Section~\ref{sec:our NN} \\
 \hline\hline
Vector representation of generators & Size of GB  & 0.270 & 0.6148\\
(data set A) & Max GB degree & Not Available & Not Available\\
\hline
Vector representation of generators & Size of GB  & 0.180 & 0.401\\
(data set B) & Max GB degree & 0.115& -0.1706\\
\hline
Features summarizing generators & Size of GB & 0.068 & 0.083\\
(data set C) & Max GB degree & 0.023 & -0.2850\\
\hline\\
\end{tabular}
\caption{$r^2$ values for various predictions for reduced reverse lexicographic Gr\"obner basis. This statistic is a numerical measure of  how close the true values lie to the predicted quantity. 
}\label{r2s}
\end{table*}

From the point of view of machine learning, both network architecture and algorithms we used in this paper are straightforward. We expect  performance improvements can be achieved  in the future, with more training and other data sets.

\section{Data and code availability}
The data we used will be made available on Zenodo shortly, in similar fashion as in \cite{MPP22}, along with the code for generating other random binomial ideal data and the code for learning and prediction. The latter is also housed on Github, under the following link: 

\url{https://sondzus.github.io/LearningGBsize/}

The GitHub page contains {\tt Jupyter} notebooks with the code for learning and prediction described in the paper.
The data files for 250,000 ideals are in the 4 text files zipped; the full data sets will be hosted on Zenodo, due to file size.
Code for generating data - random ideals and their Gr\"obner bases - is in {\tt Macaulay2}, and the site also continas small example output files.

\section*{Acknowledgements}
Shahrzad Jamshidi is with Lake Forest College, where Eric Kang was an undergraduate researcher. Sonja Petrovi\'c is with  Illinois Institute of Technology, supported by Simons Foundation's Travel Support for Mathematicians \#854770  and DoE award \#1010629. 
The authors are grateful to Despina Stasi, who helped with  the {\tt Macaulay2} code  used to generate the binomial ideal data.

\bibliographystyle{agsm}

\bibliography{MLandAlg,NNs,SparkRandomizer,SparkRandomizerWhitepaper}

\end{document}

%% file: sectionBackgroundNN.tex
With the advent of big data and the wide availability of cheap computation, there has been an explosion of work in the field of machine learning (ML), and, in particular, in the use of neural networks. Machine learning is the study of computer programs designed to ``learn,'' that is, to improve their performance through experience (e.g.\cite{Mitchell97}).
In practice the program is fed a dataset to train on (a training set), either labelled (for supervised learning) or fully unlabelled (for unsupervised learning) and must produce a labelling of the instances in the dataset. From a mathematical perspective, one can think of these  ML algorithms as producing function approximations, whether this function
describes a relationship between variables (regression) or separates the dataset (classification).

Neural networks (NNs) are supervised learning models within machine learning that draw inspiration from scientific models of biological neurons.  
They are graphically represented using directed bipartite graphs where each node represents a {\it neuron}---the entry in a vector. These neurons are organized into layers, which represent vectors. Layers which are not the first (input values) or last layer (output values) are called hidden layers.

Edges are drawn between layers to indicate which values of the previous layer are transformed to define the given neuron. We should note that these transformations are typically nonlinear. Within each transformation, there are parameters, referred to as {\it weights and biases}, that are adjusted to better fit a dataset. These adjustments are done using a form of gradient descent.

Theoretically, NNs can approximate any continuous function with compact support (as shown in \cite{HORNIK1989359}). This is a concept referred to in the ML community as ``universal approximation." In practice, they approximate a particular class of functions that are dictated by the choice of hyperparameters. These hyperparameters are often chosen and adjusted using ad-hoc, trial-and-error methods. Hyperparameters include
\begin{itemize}
  \item the network topology,
  \item the activation function,
  \item the cost function,
  \item the number of epochs,
  \item the mini-batch size,
  \item dropout thresholds, and
  \item choice of optimizer.
\end{itemize}

The {\it network topology} defines the number of layers and the number of nodes per layer. Figure \ref{neuralnetwork_layers} depicts a neural network with two hidden layers, densely connected. Each layer has three neurons. The number of neurons per layer or edges between layers can be adjusted, with the main limiting factor being computational power.

\begin{figure}
\begin{center}
\begin{tikzpicture}[scale = 2, node distance = 2cm]]
\node [draw, circle, inner sep = 2pt] (x1) {$x_1$};
\node [draw, circle, above of = x1, inner sep = 2pt] (x2) {$x_2$};
\node [draw, circle, above of = x2, inner sep = 2pt] (x3) {$x_3$};
\node [right of=x1, draw, circle, inner sep = 5pt] (node11){};
\node [right of=x2, draw, circle, inner sep = 5pt] (node12){};
\node [right of=x3, draw, circle, inner sep = 5pt] (node13){};
\node [right of=node11, draw, circle, inner sep = 5pt] (node21){};
\node [right of=node12, draw, circle, inner sep = 5pt] (node22){};
\node [right of=node13, draw, circle, inner sep = 5pt] (node23){};
\node [right of = node21, draw, circle, inner sep = 2pt] (node31) {$y_1$};
\node [right of = node22, draw, circle, inner sep = 2pt] (node32) {$y_2$};
\node [right of = node23, draw, circle, inner sep = 2pt] (node33) {$y_3$};
\foreach \i in {1, 2, 3}
{
\foreach \j in {1, 2, 3}
{
\draw[->] (x\i) -- (node1\j);
\draw[->] (node1\i) -- (node2\j);
\draw[->] (node2\i) -- (node3\j);
}
}
\node[rotate=90] at (-.5, 1) {\it Inputs};
\node[rotate=90] at (3.5, 1)(outputs) {\it Outputs};
\draw[red, thick] (.8, 2.2) --node[above]{${\bf x}_1$} (1.2, 2.2) -- (1.2, -.2) -- (.8, -.2) -- (.8, 2.2);

\draw[red, thick] (1.8, 2.2) --node[above]{${\bf x}_2$} (2.2, 2.2) -- (2.2, -.2) -- (1.8, -.2) -- (1.8, 2.2);

\end{tikzpicture}
\end{center}
\caption{Example of a (feed forward) neural network with two hidden layers. Each layer is considered dense because it is fully connected to the neighboring layers. The input layer is used to compute the values in the first hidden layer, ${\bf x}_1$. This first hidden layer is used to compute ${\bf x}_2$ and, so on. In the above diagram, the output layer, ${\bf y}$, the last hidden layer, in this case ${\bf x}_2$. See \ref{nneqs} for the mathematical encoding.}\label{neuralnetwork_layers}
\end{figure}

\begin{figure}
\begin{center}
\fbox{\begin{minipage}{.45\textwidth}
The input layer is ${\bf x} = (x_1, x_2, x_3)$ and is used to compute the first hidden layer, ${\bf x}_1 = (x_{1, 1}, x_{2, 1}, x_{3, 1})$, by 
\begin{align*}
&x_{1, 1} = \varphi({\bf w}_{1, 1} \cdot {\bf x} + b_{1, 1})\\
&x_{2, 1} = \varphi({\bf w}_{2, 1} \cdot {\bf x} + b_{2, 1})\\
&x_{3, 1} =\varphi({\bf w}_{3, 1} \cdot {\bf x} + b_{3, 1}).
\end{align*}

 The second hidden layer is defined analogously using ${\bf x}_1$ in place of ${\bf x}$.\\


The output layer consists of
\begin{align*}
&y_1 = \varphi({\bf w}_{1, 3} \cdot {\bf x}_2 + b_{1, 3})\\
&y_2 = \varphi({\bf w}_{2, 3} \cdot {\bf x}_2 + b_{2, 3})\\
&y_3 = \varphi({\bf w}_{3, 3} \cdot {\bf x}_2 + b_{3, 3})
\end{align*}
which is also the output.
\end{minipage}}
\caption{Example of a neural network with 2 hidden layers:  the  equations corresponding to the network in Figure~\ref{neuralnetwork_layers}. Hidden neuron values are denoted as $x_{i,j}$, corresponding to $i^{th}$ node of the $j^{st}$ layer. ${\bf w}_{i,j}$ and  $b_{i,j}$ are weights and biases  for the $i^{th}$ node of the $j^{th}$ layer, respectively.} \label{nneqs}
\end{center}
\end{figure}

The {\it activation function}, $\varphi$, defines the transformation represented by edges. The sum total of the edges going into a neuron represent a composition of the activation function with an affine function evaluated at the previous layer, of the form:  $$\varphi \circ g({\bf x}) = \varphi({\bf w} \cdot {\bf x} + b).$$ The variables ${\bf w}$ and $b$ are parameters referred to as the {\it weights} and {\it bias}, respectively. For Figure \ref{neuralnetwork_layers}, the diagram corresponds to the set of equations in Figure \ref{nneqs}.
The connections between layers may be {\it dense}, meaning all the values of the previous vector are used to define each entry of the next one.

Edges represent a calculation constructed using an affine function of parameters and $\varphi: \mathbb{R} \to \mathbb{R}$. Originally, models used sigmoidal functions like \begin{equation}\label{sigmoid}\varphi(x) = \frac{1}{1 + e^{-x}}\end{equation} based on the biological work by Hartline and Ratliff in \cite{hartline_ratliff_1957} and others, that were inspired by the stimulus response of a physical neuron; however, other choices of $\varphi(\cdot)$ have been explored as many are able to satisfy the requirement of universal approximation defined in \cite{cybenko, hornik}. The most widely used activation function appears to be \begin{equation} \label{relu} \varphi(x) = max\{0,x\}
\end{equation} referred to as the Rectified Linear Unit (ReLU), likely due to its speed and efficacy. \cite{Zhang2018} describes a correspondence between feed-forward neural networks with ReLU activation and tropical rational maps. They provide a meaningful discussion with respect to how these neural networks can approximate high dimensional surfaces well.

Parameters within the transformations---the weights and biases---are adjusted using an approximation of gradient descent called {\it backpropagation}. With every choice of values in the parameters, an error is measured between the estimated output (labels) and the actual output using a cost function. The gradient corresponding to the {\it cost function} is then used to adjust the parameters. An example of a cost function is the mean squared error of the difference between the estimated and actual labels.

An {\it epoch} refers to the number of times the parameters are adjusted using the entire training set. The gradient of the cost function can be measured using the entire dataset, so an epoch can be a one-step process. It is common, however, to break up the dataset into mini-batches (sometimes just called batches). Each mini-batch is used to construct a gradient and adjust the parameters. Under these conditions, an epoch is completed once all the mini-batches are used. Breaking up the training set into mini-batches reduces the memory requirements for each update of the parameters and introduces some noise in the gradient descent process; smaller mini-batches will introduce more noise. Some noise is thought to be beneficial as it discourages settling into suboptimal minima. Too much noise, however, can disrupt convergence to an optimal minimum.

Often, networks begin with dense connections between layers--meaning that every neuron in a layer is connected to all the neurons in the previous layer. A threshold can be defined, informally called {\it dropout}, that forces parameters to zero if they are below the specified threshold. Effectively, these parameters are ``dropped out" of the model thereby deleting an edge.

Optimization algorithms can be used to improve the gradient descent process in the tuning of a neural network. These work by rescaling and/or adding stochasticity to the gradient in strategic ways so as to improve convergence.